\documentclass[11pt,reqno]{amsart}

\usepackage{amssymb,upref,url,mathrsfs}
\let\mathcal\mathscr

\newtheorem{theorem}{Theorem}
\newtheorem{lemma}[theorem]{Lemma}
\newtheorem{proposition}[theorem]{Proposition}
\newtheorem{corollary}[theorem]{Corollary}

\theoremstyle{remark}

\newtheorem*{remark}{Remark}
\newtheorem*{RegularityCondition}{Regularity Condition
\textup{(\cite{HenneauxTeitelboim:QGS})}}

\newcommand{\cprime}{\/{\mathsurround=0pt$'$}}

\newcommand*{\pd}[2]{\mathchoice{\frac{\partial #1}{\partial #2}}
                    {\partial #1/\partial #2}{\partial #1/\partial #2}
                    {\partial #1/\partial #2}}

\newcommand{\eval}[2][\right]{\relax
  \ifx#1\right\relax \left.\fi#2#1\rvert}

\newcommand{\dotsa}{\cdots}

\renewcommand{\kappa}{\varkappa}
\renewcommand{\phi}{\varphi}

\newcommand{\R}{\mathbb{R}}
\newcommand{\Z}{\mathbb{Z}}

\DeclareMathOperator{\CDiff}{\mathcal{C}Dif{}f}

\DeclareFontFamily{OT1}{wncyr}{}
\DeclareFontShape{OT1}{wncyr}{m}{rm}{
   <5> <6> <7> <8> <9> gen * wncyr
   <10> <10.95> <12> <14.4> <17.28> <20.74> <24.88> wncyr10
  }{}
\DeclareSymbolFont{cyrletters}{OT1}{wncyr}{m}{rm}
\DeclareSymbolFontAlphabet{\cyrmath}{cyrletters}
\DeclareMathSymbol{\rez}{\cyrmath}{cyrletters}{"5A}

\DeclareFontFamily{OT1}{wncyi}{}
\DeclareFontShape{OT1}{wncyi}{m}{it}{
   <5> <6> <7> <8> <9> gen * wncyi
   <10> <10.95> <12> <14.4> <17.28> <20.74> <24.88> wncyi10
  }{}
\DeclareSymbolFont{cyiletters}{OT1}{wncyi}{m}{it}
\DeclareSymbolFontAlphabet{\cyimath}{cyiletters}
\DeclareMathSymbol{\Ev}{\cyimath}{cyiletters}{"03}

\providecommand{\href}[1]{}
\newcommand*{\eprint}[1]{\href{http://arXiv.org/abs/#1}%
{\begingroup \Url{arXiv:#1}}}

\begin{document}

\hfill Preprint ESI 1032

\hfill DIPS-6/2001

\hfill math.DG/0105207

\hfill To appear in Acta Appl. Math.

\bigskip

\title[Compatibility complex and the Koszul-Tate resolution]{Remarks
on two approaches to the horizontal cohomology: compatibility
complex and the Koszul-Tate resolution}
  
\author{Alexander Verbovetsky}

\address{Independent University of Moscow.  Correspondence to:\newline
A.~Verbovetsky, Profsoyuznaya~98-9-132, 117485~Moscow, Russia.}

\email{verbovet@mccme.ru}

\date{May 2001}

\thanks{Research supported in part by the GNFM of the INdAM (Italy),
the University of Lecce, the Italian MURST Young Researchers Project,
and the ESI (Vienna).}

\begin{abstract} 
  The Koszul-Tate resolution is described in the context of the
  geometry of jet spaces and differential equations.  The application
  due to Barnich, Brandt, and Henneaux of this resolution to computing
  the horizontal cohomology is analyzed.  Relations with the
  Vinogradov spectral sequence are discussed.
\end{abstract}
 
 

\maketitle

\begin{flushright}
  \begin{minipage}{.6\columnwidth}
    The Perl motto is ``There's more than one way to do it.''
    Divining how many more is left as an exercise to the reader.
    \smallskip
    \begin{flushright}
      Larry Wall, \textit{The Perl man page}
    \end{flushright}
  \end{minipage}
\end{flushright}

\section{Introduction}  

This paper is concerned with general methods for computing horizontal
(also called ``characteristic'') cohomology of systems of nonlinear
partial differential equations.  There are at least two such methods.
One stems from the fact that the horizontal cohomology is the column
$E_1^{0,\bullet}$ of the Vinogradov spectral sequence
\cite{Vinogradov:AlGFLFT,Vinogradov:SSAsNDEqAlGFLFTC,%
Vinogradov:SSLFCLLTNT} and thereby related to the terms
$E_1^{p,\bullet}$ for $p>0$.  These terms can be computed via the
compatibility complex for the linearization of the system under
consideration.  For a detailed description of this technique the
reader should consult
\cite{KrasilshchikVerbovetsky:HMEqMP,Verbovetsky:NHC}.

The second method was proposed in
\cite{BarnichBrandtHenneaux:LBCAnFIGT} and is based on the Koszul-Tate
resolution \cite{HenneauxTeitelboim:QGS}.  The purpose of the present
paper is to describe this method in the language of the geometry of
differential equations (see, e.g.,
\cite{KrasilshchikVinogradov:SCLDEqMP,Krasilshchik:GDEqCIn}) and look
into relationships between this method and the former one.

We restrict the discussion to a general theory.  As an example we
refer to the $p$-form gauge theory that was explicitly worked out in
\cite{HenneauxKnaepenSchomblond:CCFGT} by means of the Koszul-Tate
resolution and in
\cite{Verbovetsky:NHC,KrasilshchikVerbovetsky:HMEqMP} by means of the
Vinogradov spectral sequence.

\section{On compatibility operators}
\label{sec:comp-oper}

We assume that the reader is familiar with the geometry of jet spaces
and differential equations, including the horizontal cohomology, the
Vinogradov spectral sequence, the compatibility complex, and the
$k$-line theorem.  This material can be found in
\cite{KrasilshchikVerbovetsky:HMEqMP,Verbovetsky:NHC}.  Throughout the
paper we use definitions, notation, and terminology from these works.

Let $\pi\colon E\to M$ be a vector bundle and $\pi_\infty\colon
J^\infty(\pi)\to M$ the associated infinite jet bundle.  The standard
coordinates on the space~$J^\infty(\pi)$ are coordinates $x_i$ of the
manifold~$M$ and infinitely many coordinates $u^j_\sigma$,
corresponding to partial derivatives of sections $u=(u^1,\dotsc,u^m)$
of the bundle~$\pi$.
\begin{remark}
  All constructions of the present paper can be readily generalized to
  the case of vector bundle~$\pi$ with super fibers.
\end{remark}

Let $\alpha_1$ be a vector bundle over~$M$,
$P_1=\Gamma(\pi_\infty^*(\alpha_1))$ a horizontal module on
$J^\infty(\pi)$ (as usual, a module is always the module of sections
of a vector bundle; $\Gamma$ denotes the functor that takes bundles to
the corresponding modules), $\mathcal{E}=\{F=0\}\subset J^k(\pi)$ a
formally integrable equation, determined by a section $F\in P_1$, and
$\mathcal{E}^\infty\subset J^\infty(\pi)$ its infinite prolongation.
\begin{remark}
  Physicists called coordinates along fibers of~$\alpha_1$
  \emph{antifields} and say that they have the \emph{antighost number}
  $1$ (hence the notation~$P_1$).
\end{remark}

\begin{RegularityCondition}
  We shall assume that there exists an open submanifold
  $\mathcal{U}\subset J^\infty(\pi)$ and an isomorphism
  $\nu\colon\mathcal{U}\to\mathcal{E}^\infty\times\mathcal{V}$, where
  $\mathcal{V}$ is a star-shaped neighborhood of the zero in
  $\R^\infty$, such that $\mathcal{U}\supset\mathcal{E}^\infty$,
  $\nu(\theta)=(\theta,0)$ for any point
  $\theta\in\mathcal{E}^\infty$, and the composition
  $v\colon\mathcal{U}\xrightarrow{\nu}
  \mathcal{E}^\infty\times\mathcal{V}\to\mathcal{V}$ of the
  isomorphism~$\nu$ and the projection on the second factor has the
  form
  \begin{displaymath}
    v=(v_1,v_2,v_3,\dotsc,v_s,\dotsc),
  \end{displaymath}
  where $v_s=D_{\sigma_s}(F_{i_s})$.
\end{RegularityCondition}

Because of the regularity condition, a function~$f$ on~$\mathcal{U}$
vanishes on~$\mathcal{E}^\infty$ if and only if it has the form
$f=\Delta(F)$ for some operator $\Delta\in\CDiff(P_1,\mathcal{F})$,
where $\mathcal{F}$ is the algebra of functions on~$\mathcal{U}$.

From here on, we fix a manifold~$\mathcal{U}$ and will consider it
instead of the whole jet space $J^\infty(\pi)$.

\begin{lemma}
  Let $Q=\Gamma(\pi_\infty^*(\alpha))$ be a horizontal
  module\textup{;} then the kernel of the map
  \begin{displaymath}
    Q\xrightarrow{}\CDiff(\eval{\kappa}_{\mathcal{E}^\infty},
    \eval{Q}_{\mathcal{E}^\infty}),\qquad
    q\mapsto\eval{\ell_q}_{\mathcal{E}^\infty},
  \end{displaymath}
  where $\kappa=\Gamma(\pi_\infty^*(\pi))$ and $\ell_q$ is the
  linearization of~$q$\textup{,} has the form $\mathcal{I}^2Q\oplus
  Q'\subset Q$, where $\mathcal{I}=\{\,\Delta(F)
  \mid\CDiff(P_1,\mathcal{F})\,\} \subset\mathcal{F}$ is the ideal of
  the equation~$\mathcal{E}^\infty$ and $Q'=\{\,\pi_\infty^*(s)\mid
  s\in\Gamma(\alpha)\,\}$\textup{.}  \textup{(}In other
  words\textup{,} the set $\mathcal{I}^2Q\subset Q$ consists of
  elements of the form $\Delta(F,F)$ for
  $\Delta\in\CDiff(P_1,\CDiff(P_1,Q))$\textup{.)}
\end{lemma}
\begin{proof}
  It suffices to consider the case $Q=\mathcal{F}$ and to prove the
  statement in a local chart.  Take a function~$f\in\mathcal{F}$ such
  that $\ell_f=0$.  Then $\pd{f}{u_\sigma^j}\in\mathcal{I}$.  Now let
  us change coordinates to $(x_i,v_s,w_t)$, where $v_s$ are the
  coordinates along~$\mathcal{V}$ from the regularity condition and
  $w_t$ are arbitrary coordinates along~$\mathcal{E}^\infty$.  We have
  $\pd{f}{v_s}\in\mathcal{I}$ and $\pd{f}{w_t}\in\mathcal{I}$.  Write
  the function~$f$ in the form
  \begin{displaymath}
    f=f_0(x_i,w_t)+\sum_sf_s(x_i,w_t)v_s+f',
  \end{displaymath}
  where $f'\in\mathcal{I}^2$.  Then
  $\pd{f}{v_s}=f_s+\pd{f'}{v_s}\in\mathcal{I}$, hence $f_s=0$, and
  $\pd{f}{w_t}=\pd{f_0}{w_t}+\pd{f'}{w_t}\in\mathcal{I}$, so that
  $\pd{f}{w_t}=0$.
  Thus, $f\in\mathcal{I}^2\oplus\pi_\infty^*(C^\infty(M))$.
\end{proof}

\begin{proposition}
  An operator $\Delta\in\CDiff(\eval{P_1}_{\mathcal{E}^\infty},
  \eval{P_2}_{\mathcal{E}^\infty})$ is the compatibility operator for
  $\eval{\ell_{F}}_{\mathcal{E}^\infty}$ if and only if for each
  operator $\nabla\colon P_1\to Q$ such that $\nabla(F)=0$ we have
  $\eval{\nabla}_{\mathcal{E}^\infty}=\square\circ\Delta$ for an
  operator $\square\in\CDiff(\eval{P_2}_{\mathcal{E}^\infty},
  \eval{Q}_{\mathcal{E}^\infty})$\textup{.}
\end{proposition}
\begin{proof}
  Suppose that $\Delta\in\CDiff(\eval{P_1}_{\mathcal{E}^\infty},
  \eval{P_2}_{\mathcal{E}^\infty})$ is the compatibility operator for
  $\eval{\ell_{F}}_{\mathcal{E}^\infty}$. Linearizing the equality
  $\nabla(F)=0$, we get $\eval{\nabla}_{\mathcal{E}^\infty}
  \circ\eval{\ell_{F}}_{\mathcal{E}^\infty}=0$.  Hence,
  $\eval{\nabla}_{\mathcal{E}^\infty}=\square\circ\Delta$.
  
  Conversely, consider an operator $\nabla\in\CDiff(P_1,Q)$
  satisfying the condition $\eval{\nabla}_{\mathcal{E}^\infty}\circ
  \eval{\ell_F}_{\mathcal{E}^\infty}=0$.  Since
  $\eval{\ell_{\nabla(F)}}_{\mathcal{E}^\infty}
  =\eval{\nabla}_{\mathcal{E}^\infty}\circ
  \eval{\ell_F}_{\mathcal{E}^\infty}=0$, we see that
  $\nabla(F)=\nabla'(F,F)$.  Define the operator
  $\nabla_1\in\CDiff(P_1,Q)$ by the equality
  $\nabla_1(p)=\nabla(p)-\nabla'(F,p)$.  We have $\nabla_1(F)=0$, so
  that $\eval{\nabla}_{\mathcal{E}^\infty}
  =\eval{\nabla_1}_{\mathcal{E}^\infty}=\square\circ\Delta$.  This
  completes the proof.
\end{proof}

An equation is called \emph{normal} if the compatibility operator for
the operator $\eval{\ell_F}_{\mathcal{E}^\infty}$ is trivial.
\begin{corollary}
  The equation $\mathcal{E}^\infty$ is normal if and only if for each
  operator $\nabla\colon P_1\to Q$ the equality $\nabla(F)=0$ implies
  $\eval{\nabla}_{\mathcal{E}^\infty}=0$\textup{.}
\end{corollary}
\begin{corollary}[from the proof]
  \label{sec:comp-compl-3}
  Suppose that $\Delta\in\CDiff(\eval{P_1}_{\mathcal{E}^\infty},
  \eval{P_2}_{\mathcal{E}^\infty})$ is the compatibility operator for
  $\eval{\ell_F}_{\mathcal{E}^\infty}$\textup{;} then there exist an
  extension $\tilde{\Delta}\in\CDiff(P_1,P_2)$ of~$\Delta$ such that
  $\tilde{\Delta}(F)=0$\textup{.}
\end{corollary}

\section{The Koszul-Tate resolution}
\label{sec:kosz-tate-resol}

Due to the last Corollary there is an operator $\Delta_1\colon P_1\to
P_2$ such that $\eval{\Delta_1}_{\mathcal{E}^\infty}$ is the
compatibility operator for $\eval{\ell_F}_{\mathcal{E}^\infty}$ and
$\Delta_1(F)=0$.  Pick up such an operator and construct the
compatibility complex for it:
\begin{displaymath}
  P_1\xrightarrow{\Delta_1}P_2\xrightarrow{\Delta_2}\dotsa
  \xrightarrow{\Delta_{k-2}}P_{k-1}\xrightarrow{}0,
\end{displaymath}
where $P_i=\Gamma(\alpha_i)$ for some vector bundles $\alpha_i\colon
V_i\to \mathcal{U}$.
\begin{remark}
  Such a complex may not exist (modules $P_i$ may not be projective).
  In \cite{HenneauxTeitelboim:QGS} the assumption that it exists is
  termed as ``off-shell reducibility.''  The situation then this
  complex exists only on~$\mathcal{E}^\infty$ is called ``on-shell
  reducibility.''  So, we here require the off-shell reducibility.
\end{remark}
Take the direct sum
$\alpha=\bigoplus_{i\geq1}(\alpha_{2i-1}^\Pi\oplus\alpha_{2i})$, where
$\alpha_{2i-1}^\Pi$ means the bundle $\alpha_{2i-1}$ with reversed
parity of fibers, and consider the supermanifold of horizontal jets of
its sections $\bar{J}^\infty(\alpha)$.

Each Cartan field on~$\mathcal{U}$ can be naturally lifted to
$\bar{J}^\infty(\alpha)$; these liftings span the Cartan distribution
on the horizontal jets~$\bar{J}^\infty(\alpha)$.  It is not hard to
check that all the theory of jet spaces can be carried over to the
space~$\bar{J}^\infty(\alpha)$.  In coordinate language, the
horizontal jet space is a jet space equipped with extra base
coordinates~$u_\sigma^j$ (so that all functions depend on them as on
parameters), with the total derivatives
$D_i=\pd{}{x_i}+\sum_{j,\sigma}u_{\sigma i}^j\pd{}{u_{\sigma}^j}$ in
place of the partial derivatives $\pd{}{x_i}$.

Now let us pull the element~$F$ and all operators~$\Delta_i$ back
to~$\bar{J}^\infty(\alpha)$.  We shall treat them as elements of the
module $\kappa(\alpha)=\Gamma(\alpha^*(\alpha))$.  Consider the
element $\Phi=F+\Delta_1+\dotsb+\Delta_{k-2}\in\kappa(\alpha)$.
\begin{proposition}
  The odd evolutionary vector field $\delta=\Ev_\Phi$ is a
  differential\textup{:} $\delta^2=0$\textup{.}
\end{proposition}
\begin{proof}
  Since $\delta$ is a vector field and
  $\eval{\delta}_{\mathcal{F}}=0$, it is sufficient to
  evaluate~$\delta^2$ on functions linear along the fibers of the
  natural projection~$\alpha_\infty\colon\bar{J}^\infty(\alpha)\to
  \mathcal{U}$.  Such functions can be naturally identified with
  $\mathcal{C}$-differential operators belonging to
  $\CDiff(P,\mathcal{F})$, where $P=\Gamma(\alpha)
  =\bigoplus_{i\geq1}(\Gamma(\alpha_{2i-1}^\Pi)\oplus
  \Gamma(\alpha_{2i}))$ is a graded $\mathcal{F}$-module.  Define an
  odd $\mathcal{C}$-differential operator $\Delta\colon P\to P$ by the
  formula
  \begin{displaymath}
    \Delta(p_1,\dotsc,p_{k-2},p_{k-1})
    =(0,\Delta_1(p_1),\dotsc,\Delta_{k-2}(p_{k-2})),
  \end{displaymath}
  where $p_i\in P_i=\Gamma(\alpha_i)$.  Trivially, $\Delta^2=0$ and
  $\Delta(F)=0$.  It is easy to see that if
  $\nabla\in\CDiff(P,\mathcal{F})$ then
  $\delta(\nabla)=\nabla(F)+\nabla\circ\Delta$, so that
  $\delta^2(\nabla)=\nabla\circ\Delta^2=0$.
\end{proof}

\begin{remark}
  Elements of a module over a jet space can be naturally thought of as
  graded symmetric multilinear operators.  We used such an
  identification (applied to linear functions) in the last proof. This
  is a simple but very useful trick.  To facilitate further
  observations, it might be well always to bear in mind this
  identification.
\end{remark}

Denote by $\mathcal{F}^\mathrm{pol}(\alpha)$ the subalgebra of the
algebra~$\mathcal{F}(\alpha)$ of functions on $\bar{J}^\infty(\alpha)$
that consists of functions polynomial along the fibers of the
projection~$\alpha_\infty\colon\bar{J}^\infty(\alpha)\to \mathcal{U}$.
We supply the algebra $\mathcal{F}^\mathrm{pol}(\alpha)$ with a
$\Z$-grading, $\mathcal{F}^\mathrm{pol}(\alpha)
=\bigoplus_{i\geq0}\mathcal{F}^\mathrm{pol}_i(\alpha)$, called the
\emph{antighost number}, such that fiberwise linear functions on
$\bar{J}^\infty(\alpha_i)$ have antighost number~$i$ and functions on
$\mathcal{U}$ have antighost number zero.  (Thus the parity of a
function is equal to its antighost number modulo~$2$.)  The
differential~$\delta$ reduces the antighost number by~$1$ and so we
have the complex
\begin{displaymath}
  0\xleftarrow{}\mathcal{F}\xleftarrow{\delta}
  \mathcal{F}^\mathrm{pol}_1(\alpha)\xleftarrow{\delta}
  \mathcal{F}^\mathrm{pol}_2(\alpha)\xleftarrow{\delta}\dotsa
\end{displaymath}
called \emph{the Koszul-Tate complex} \cite{HenneauxTeitelboim:QGS}.
\begin{theorem}
  \label{sec:kosz-tate-resol-1}
  The Koszul-Tate complex is a resolution\textup{,} with the zero
  homology isomorphic to the algebra~$\mathcal{F}(\mathcal{E})$ of
  functions on the equations~$\mathcal{E}^\infty$\textup{.}  The
  homology of the differential group $(\mathcal{F}(\alpha),\delta)$ is
  equal to~$\mathcal{F}(\mathcal{E})$ as well\textup{.}
\end{theorem}
\begin{proof}
  See~\cite{HenneauxTeitelboim:QGS}.
\end{proof}

Let us consider the bicomplex
$(\mathcal{F}^\mathrm{pol}_p(\alpha)\otimes\bar{\Lambda}^q,
\delta,(-1)^q\bar{d})$, where $\bar{\Lambda}^q$ is the module of
horizontal $q$-forms on~$\mathcal{U}$.  The horizontal
differential~$\bar{d}$ is well-defined, since there is a natural
inclusion $\mathcal{F}^\mathrm{pol}_p(\alpha)
\otimes\bar{\Lambda}^q\subset\bar{\Lambda}^q(\alpha)$, where
$\bar{\Lambda}^q(\alpha)$ is the module of horizontal $q$-forms
on~$\bar{J}^\infty(\alpha)$.  Using this bicomplex and
Theorem~\ref{sec:kosz-tate-resol-1}, we immediately conclude that for
$0\leq q\leq n-1$ the group
$\bar{H}^q(\mathcal{E})/\bar{H}^q(\mathcal{U})$ is isomorphic to the
$(n-q)^{\mathrm{th}}$ homology group of the complex
\begin{equation}
  \label{eq:1}
  \bar{H}^n_0(\alpha)\xleftarrow{\delta}
  \bar{H}^n_1(\alpha)\xleftarrow{\delta}
  \bar{H}^n_2(\alpha)\xleftarrow{\delta}\dotsa,
\end{equation}
where $\bar{H}^n_i(\alpha)
=H^n(\mathcal{F}^\mathrm{pol}_i(\alpha)\otimes
\bar{\Lambda}^q,\bar{d})$.  Application of the one-line theorem for
$\mathcal{U}$ yields the following reformulation of this result.
\begin{theorem}
  \label{sec:kosz-tate-resol-3}
  For $0\leq q\leq n-1$ there is the isomorphism
\begin{displaymath}
  \bar{H}^q(\mathcal{E})/H^q(M)
  =H_{n-q}(\bar{H}^n_\bullet(\alpha),\delta).
\end{displaymath}
\end{theorem}

Our main concern now is to compute the homology of
complex~\eqref{eq:1}.  To this end it is convenient to drop the
zero-term $\bar{H}^n_0(\alpha)$ and to embed the rest of the complex
to the complex $(\hat{\kappa}(\alpha),\delta+\ell_\Phi^*)$ by means of
the Euler operator
\begin{equation}
  \label{eq:2}
  \bigoplus_{i\geq1}\bar{H}^n_i(\alpha)\to\hat{\kappa}(\alpha)
\end{equation}
on $\bar{J}^\infty(\alpha)$.
\begin{remark}
  We use here two facts:
  \begin{enumerate}
  \item the map $\delta+\ell_\Phi^*\colon\hat{\kappa}(\alpha)\to
    \hat{\kappa}(\alpha)$ is a differential:
    $(\delta+\ell_\Phi^*)^2=0$;
  \item the Euler operator~\eqref{eq:2} is a chain map.
  \end{enumerate}
  Both of these statements follow from the infinitesimal Stokes
  formula in the term~$E_1$ of the Vinogradov spectral sequence
  on~$\bar{J}^\infty(\alpha)$
  (see~\cite[Sec.~9.10]{Vinogradov:SSLFCLLTNT}\footnote{Note a
  misprint in~\cite[eq.~$(9.10.5)$]{Vinogradov:SSLFCLLTNT}.  The
  correct formula (in the notation of~\cite{Vinogradov:SSLFCLLTNT}) is
  \begin{displaymath}
    \chi\{\phi\}
    =l_\phi(\chi)+l_\chi^*(\phi)=\rez_\chi(\phi)+l_\chi^*(\phi).
  \end{displaymath}
}).
\end{remark}
\begin{proposition}
  The operator
  $\ell_\Phi^*\colon\hat{\kappa}(\alpha)\to\hat{\kappa}(\alpha)$ is a
  differential\textup{,} i\textup{.}e\textup{.,}\linebreak
  $(\ell_\Phi^*)^2=0$\textup{.}
\end{proposition}
\begin{proof}
  Pick up an element $\theta\in\hat{\kappa}(\alpha)$.  It can be
  thought of as a nonlinear operator $\theta\colon P\to\hat{P}$.
  Obviously,
  $\ell_\Phi^*(\theta)=(-1)^{p(\theta)}\Delta^*\circ\theta$, where
  $p(\theta)$ is the parity of~$\theta$.  Hence,
  $(\ell_\Phi^*)^2(\theta) =-(\Delta^*)^2\circ\theta=0$.
\end{proof}
\begin{remark}
  Note that the operator $\Delta^*\colon\hat{P}\to\hat{P}$ has the
  form
  \begin{displaymath}
    \Delta^*(\hat{p}_1,\hat{p}_2,\dotsc,\hat{p}_{k-1})
    =(-\Delta^*_1(\hat{p}_2),\Delta^*_2(\hat{p}_3),\dotsc,
    (-1)^{k-2}\Delta^*_{k-2}(\hat{p}_{k-1})),
  \end{displaymath}
  where $\hat{p}_i\in\hat{P}_i$.
\end{remark}

Thus, we see that the differentials $\delta$ and $\ell_\Phi^*$ define
a bicomplex structure in~$\hat{\kappa}(\alpha)$.  The same holds true
for the polynomial part $\hat{\kappa}^{\mathrm{pol}}(\alpha)
=\mathcal{F}^{\mathrm{pol}}(\alpha)\otimes\hat{P}$. This enables an
easy computation of the homology of the total complex
$(\hat{\kappa}(\alpha),\delta+\ell_\Phi^*)$ and its polynomial
subcomplex.
\begin{proposition}
  The homologies of complexes
  \begin{displaymath}
    (\hat{\kappa}(\alpha),\delta+\ell_\Phi^*)\qquad\text{and}\qquad
    (\hat{\kappa}^{\mathrm{pol}}(\alpha),\delta+\ell_\Phi^*)    
  \end{displaymath}
  coincide and equal the homology of the complex
  \begin{equation}
    \label{eq:3}
    0\xleftarrow{}\hat{P}_1\xleftarrow{\Delta_1^*}
    \hat{P}_2\xleftarrow{\Delta_2^*}\dotsa
    \xleftarrow{\Delta_{k-2}^*}\hat{P}_{k-1}\xleftarrow{}0\qquad\text{on
    the equation~$\mathcal{E}^\infty$}.
  \end{equation}
\end{proposition}
\begin{corollary}
  \label{sec:kosz-tate-resol-2}
  Each cycle $\psi$ belonging to $\hat{\kappa}^\mathrm{pol}_+(\alpha)
  =\bigl(\bigoplus_{i\geq1}^\infty\mathcal{F}_i^{\mathrm{pol}}
  (\alpha)\bigr)\otimes\hat{P}$ is a boundary\textup{.}
\end{corollary}

From the one-line theorem for $\bar{J}^\infty(\alpha)$ it follows that
the image of the mapping~\eqref{eq:2} coincides with the complex
$(\hat{\kappa}^{\mathrm{sp}}(\alpha),\delta+\ell_\Phi^*)$, where
\begin{displaymath}
  \hat{\kappa}^{\mathrm{sp}}(\alpha)=
  \{\,\psi\in\hat{\kappa}^{\mathrm{pol}}(\alpha)
  \mid\ell_\psi^*=\ell_\psi\,\}.
\end{displaymath}
\begin{remark}
  It should be noted that $\hat{\kappa}^{\mathrm{sp}}(\alpha)$ does
  not inherit the bicomplex structure of
  $\hat{\kappa}^{\mathrm{pol}}(\alpha)$.
\end{remark}

The space~$\hat{\kappa}^{\mathrm{sp}}(\alpha)$ can be expanded in the
sum
\begin{displaymath}
  \hat{\kappa}^{\mathrm{sp}}(\alpha)
  =\hat{P}\oplus\hat{\kappa}^{\mathrm{sp}}_+(\alpha),
\end{displaymath}
where $\hat{\kappa}^{\mathrm{sp}}_+(\alpha)=\hat{\kappa}^{\mathrm{sp}}
(\alpha)\cap\hat{\kappa}^{\mathrm{pol}}_+(\alpha)$.
\begin{proposition}
  Each cycle $\psi$ belonging to $\hat{\kappa}^\mathrm{sp}_+(\alpha)$
  is a boundary in the complex
  $(\hat{\kappa}^{\mathrm{sp}}(\alpha),\delta+\ell_\Phi^*)$\textup{.}
\end{proposition}
\begin{proof}
  We prove the statement by identifying in a natural way elements
  $\psi\in\hat{\kappa}^{\mathrm{pol}}(\alpha)$ with graded symmetric
  multilinear $\mathcal{C}$-differential operators
  \begin{displaymath}
    \nabla_\psi\colon P\times\dotsb\times P\to\hat{P}.
  \end{displaymath}
  It is not hard to check that
  $\psi\in\hat{\kappa}^{\mathrm{sp}}(\alpha)$ if and only if the
  corresponding operator~$\nabla_\psi$ is selfadjoint:
  $\nabla_\psi=\nabla_\psi^*$ (since $\nabla_\psi$ is symmetric,
  self-adjointness in one argument implies self-adjointness in the
  other arguments).  Consider the projector
  $S\colon\hat{\kappa}^{\mathrm{pol}}(\alpha)\to
  \hat{\kappa}^{\mathrm{sp}}(\alpha)$ given by
  $S(\nabla)=(\nabla+\nabla^*)/2$.  Obviously, $S^2=S$ and
  $S\circ(\delta+\ell_\Phi^*)=(\delta+\ell_\Phi^*)\circ S$.  Now, if a
  cycle $\nabla$ belongs to $\hat{\kappa}^{\mathrm{sp}}_+(\alpha)$
  then by Corollary~\ref{sec:kosz-tate-resol-2} we have
  $\nabla=(\delta+\ell_\Phi^*)\nabla'$.  This gives
  $\nabla=S(\nabla)=S(\delta+\ell_\Phi^*)\nabla'
  =(\delta+\ell_\Phi^*)(S(\nabla'))$, which is the desired conclusion.
\end{proof}

\begin{corollary}
  Complex~\eqref{eq:1} is exact in terms $\bar{H}^n_i(\alpha)$ for
  $i\geq k$\textup{.}
\end{corollary}
\begin{corollary}
  If $0\leq q\leq n-k$ then $\bar{H}^q(\mathcal{E})=H^q(M)$\textup{.}
\end{corollary}

Take a cycle $\psi\in\hat{\kappa}^{\mathrm{sp}}(\alpha)$ and write it
in the form $\psi=\psi_0+\psi_1$, where $\psi_0\in\hat{P}$ and
$\psi_1\in\hat{\kappa}^{\mathrm{sp}}_+(\alpha)$.  Obviously,
$\eval{\psi_0}_{\mathcal{E}^\infty}$ is a cycle of
complex~\eqref{eq:3}.  So, we have constructed a map from homology of
complex~\eqref{eq:1} to homology of complex~\eqref{eq:3}.

Let us examine this map more fully in dimension~$1$.  We have
$\bar{H}^n_1(\alpha)=\hat{P}_1$ and the map being considered is just
the restriction $\hat{P}_1\to\eval{\hat{P}_1}_{\mathcal{E}^\infty}$.
Further, an element $\beta\in\bar{H}^n_1(\alpha)=\hat{P}_1$ is a cycle
if and only if the horizontal cohomology class of the $n$-form
$\langle\beta,F\rangle$ is trivial, so that
$0=\ell_{\langle\beta,F\rangle}^*(1)=\ell^*_\beta(F)+\ell^*_F(\beta)$
and, hence, $\eval{\ell^*_F}_{\mathcal{E}^\infty}
(\eval{\beta}_{\mathcal{E}^\infty})=0$.  Thus, we get a map from
homology of complex~\eqref{eq:1} to homology of the complex
\begin{equation}
  \label{eq:6}
  \kappa\xleftarrow{\ell^*_F}\hat{P}_1
  \xleftarrow{\Delta_1^*}\hat{P}_2\xleftarrow{\Delta_2^*}\dotsa
  \xleftarrow{\Delta_{k-2}^*}\hat{P}_{k-1}\xleftarrow{}0\qquad\text{on
  the equation~$\mathcal{E}^\infty$}.
\end{equation}

The application of Theorem~\ref{sec:kosz-tate-resol-3} yields the map
\begin{equation}
  \label{eq:4}
  \bar{H}^q(\mathcal{E})\to H_{n-q}(\hat{P}_\bullet,\Delta^*_\bullet).
\end{equation}
It is straightforward to check that this map coincides with the
differential $d_1^{0,q}\colon E_1^{0,q}=\bar{H}^q(\mathcal{E})\to
E_1^{1,q}=H_{n-q}(\hat{P}_\bullet,\Delta^*_\bullet)$ of the Vinogradov
spectral sequence on~$\mathcal{E}^\infty$.

\section{A comparison}
\label{sec:comparison}

As it is seen from the last paragraphs of the previous section, the
homology of complex~\eqref{eq:6} is essential to computing the
horizontal cohomology via the Koszul-Tate resolution, similarly to
what happens when using the Vinogradov spectral sequence.  This
bridges the gap between two approaches.  Otherwise they are diverged
considerably.

As an example, let us discuss the computation of conservation laws for
a normal equation (i.e., such that $\Delta_1=0$).  The application of
the Koszul-Tate resolution gives the following result.
\begin{proposition}
  \label{sec:comparison-1}
  Let $\mathcal{E}^\infty$ be a normal equation\textup{.}  The space
  of conservation laws of~$\mathcal{E}^\infty$ is isomorphic to the
  space of solutions of the equation
  \begin{displaymath}
    \ell_F^*(\psi)+\ell_\psi^*(F)=0,
  \end{displaymath}
  where $\psi\in\hat{P}_1/\Theta$ and
  $\Theta=\{\,\square(F)\in\hat{P}_1\mid
  \square\in\CDiff(P_1,\hat{P}_1),\ \square=-\square^*\,\}$\textup{.}
\end{proposition}
Another result is obtained by means of the Vinogradov spectral
sequence.
\begin{proposition}
  \label{sec:comparison-2}
  Let $\mathcal{E}^\infty$ be a normal equation\textup{.}  The space
  of conservation laws of~$\mathcal{E}^\infty$ is a subset of the
  space of solutions of the equation
  \begin{equation}
    \label{eq:5}
    \ell_F^*(\psi)=0 \qquad{\text{on~$\mathcal{E}^\infty$,}}
  \end{equation}
  where
  $\psi\in\eval[\bigr]{\hat{P}_1}_{\mathcal{E}^\infty}$\textup{.}  A
  solution~$\psi$ corresponds to a conservation law if and only if
  on~$\mathcal{E}^\infty$ there exists a selfadjoint
  $\mathcal{C}$-differential operator $\nabla\colon P_1\to\hat{P}_1$
  such that
  \begin{displaymath}
    \ell_\psi+\Delta^*=\nabla\circ\ell_F,\quad\nabla^*=\nabla
    \qquad{\text{on~$\mathcal{E}^\infty$,}}
  \end{displaymath}
  where $\Delta\colon P_1\to\hat{\kappa}$ is a
  $\mathcal{C}$-differential operator satisfying on~$\mathcal{U}$
  the equality $\ell_F^*(\psi)=\Delta(F)$\textup{.}
\end{proposition}

Thus, both of these Propositions say that to compute conservation laws
of a normal equation we should start with solving
equation~\eqref{eq:5}.  Proposition~\ref{sec:comparison-2} also says
that if $\psi$ vanishes on~$\mathcal{E}^\infty$ then the corresponding
conservation law is trivial (essentially, this is the basic content of
the Proposition).  Proposition~\ref{sec:comparison-1} implies a weaker
result: it guaranteers triviality of conservation laws that correspond
to elements~$\psi$ of the form $\psi=\square(F)$ for skew-adjoint
operators~$\square$ only.

\section*{Acknowledgements}

The author wishes to extend his thanks to F.~Brandt,
I.~S.~Krasil{\cprime}shchik, and A.~M.~Vinogradov for many helpful
discussions.  He is also grateful to the Erwin Schr{\"o}dinger
Institute in Vienna, where a part of this research was completed, for
kind hospitality.

It is a pleasure to thank J.~Stasheff for his interest to the text and
useful suggestions, which were incorporated in the second version of
this eprint.



\end{document}